\documentclass[10pt,conference,twocolumn]{IEEEtran}
\usepackage{amsfonts}
\usepackage{bbm}
\usepackage{tipa}
\usepackage{mathrsfs}
\usepackage{amsfonts,amsmath,amsthm,latexsym,amssymb}
\usepackage{graphicx}
\usepackage{subfigure}
\usepackage{stfloats}
\usepackage{url}

\begin{document}

\title{Electricity Cost Minimization for a Residential Smart Grid with Distributed Generation and Bidirectional Power Transactions}
\renewcommand{\baselinestretch}{1}
\author{
\authorblockN{Yi~Liu}
\authorblockA{ Singapore University of \\ Technology and Design,
Singapore
\\Email: yi\_liu@sutd.edu.sg} \and
\authorblockN{Naveed~Ul~Hassan}
\authorblockA{Department of Electrical Engineering,\\ SSE, LUMS, Lahore
\\Email:naveed.hassan@yahoo.com} \and

\authorblockN{Shisheng~Huang, Chau~Yuen}
\authorblockA{ Singapore University of \\ Technology and Design,
Singapore
\\Email:\{shisheng\_huang,yuenchau\}@sutd.edu.sg}

 }\maketitle

\begin{abstract}
In this paper, we consider the electricity cost minimization problem
in a residential network where each community is equipped with a
distributed power generation source and every household in the
community has a set of essential and shiftable power demands. We
allow bi-directional power transactions and assume a two-tier
pricing model for the buying and selling of electricity from the
grid. In this situation, in order to reduce the cost of electricity
we are required to make, 1) Scheduling decisions for the shiftable
demands, 2) The decisions on the amount of energy purchased from the
gird by the users, 3) The decisions on the amount of energy sold to
the grid by the users. We formulate a global centralized
optimization problem and obtain the optimal amount of electricity
consumed, sold and purchased for each household, respectively by
assuming the availability of all current and future values of
time-varying parameters. In reality, the lack of perfect information
hampers the implementation of such global centralized optimization.
Hence, we propose a distributed online algorithm which only requires
the current values of the time-varying supply and demand processes.
We then compare and determine the tradeoff between both
formulations. Simulation results show that the proposed schemes can
provide effective management for household electricity usage.
\end{abstract}

\section{Introduction}

\IEEEPARstart {W}hile the traditional electricity grid has
maintained its current state for many decades, there is a growing
need to transit to a newer state with increasing demands from
renewable energy integration and reducing dependency on fossil fuel.
National governments and relevant industries have been making
significant efforts in the development of future electrical grids or
smart grids. A smart grid \cite{SG1}-\cite{SG3} is known as an
intelligent electricity network that uses two-way flows of
electricity and information to create an automated and distributed
power delivery system which can help save energy, reduce costs,
increase reliability and transparency.

There are two significant advantages of the smart grid: improved
integration for renewable energy (e.g. wind, solar energy) and
increased capability for demand-side management. In residential
households, the electricity generated from renewable distributed
plants such as solar panels can help the consumer save on
electricity costs and in systems supporting sell back policies like
net metering, residents can even turn a profit by selling the extra
renewable energy back to the power grid \cite{sell back}. While
typical renewable energy sources offer a cheaper and cleaner energy
supply, they introduce supply uncertainty due to the volatility of
their generation \cite{mixed energy}. It is therefore important to
include a complementary technology \cite{shisheng}, or provide a
mixture of traditional power to mitigate or cancel the volatility
\cite{renewable1} \cite{renewable2}.

Another way to mitigate this problem is to manage this supply
variability with demand side flexibility \cite{Neely1} \cite{Neely2}.
In \cite{Neely1}, the authors suggested that consumers who could
tolerate serving delays within a pre-agreed deadline can be offered electricity at a lower
price by renewable power suppliers. This essentially allows a
lower price of renewable energy to consumers willing to provide this
extra time flexibility. The renewable power supplier can now use
this flexibility to deliver the energy when it is available. Based
on this model, the authors in \cite{Neely2} investigated the problem
of allocating energy from renewable sources to flexible consumers in
electricity markets. The renewable plant must serve consumers within
a specified delay window, and incurs a cost of drawing energy from
other (possibly non-renewable) sources if its own supply is not
sufficient to meet the deadline. An natural extension to this work involving
consumer flexibility is the consideration of a distributed renewable energy source
that allows consumers to not only consume power from the grid but sell it back too.

In this paper, we consider the electricity cost minimization problem
with a renewable based Distributed Generation (DG) and bidirectional power transactions in
a residential smart grid. Each household is equipped with a mixed
renewable power source which consists of both wind and solar generators.
There are two types of appliances in each home: essential
and shiftable appliances, which leads to essential and shiftable
power demands in the entire residential system. We allow
bi-directional power transactions and assume a dynamic pricing model
for the buying and selling of electricity from the grid. In order to
reduce the cost of electricity we are required to make, 1)
Scheduling decisions for the shiftable demands, 2) The decisions on
the amount of energy purchased from the grid by the users, 3) The
decisions on the amount of energy sold to the grid by the users.

We first formulate a global optimization problem and obtain the
optimal amount of electricity consumed, sold and purchased for each
household, by assuming the availability of all current
and future values of time-varying parameters. Consequently, we
propose a distributed online algorithm which only requires the
current values of the time-varying supply and demand processes.
A comparison of the simulation results show that the proposed schemes can provide
effective management for household electricity usage.

The rest of the paper is organized as follows. In Section II, the
system models are introduced. The formulation for the global optimization
problem with perfect information is proposed in Section III.
The  distributed online algorithm with Lyapunov optimization is presented in Section IV. Section V
presents the simulation results of the proposed schemes. Finally, we
conclude the paper in Section VI.

\section{System Model}
\subsection{Mixed Distributed Renewable Energy Source}
In this subsection, we present the models for the distributed energy
sources modelled as a combination of both a  photovoltaic (PV)
module and a wind turbine, modelled separately$^*$.

\subsubsection{Output Power of a PV Module}
The irradiance data for a particular day in each season can be
described as a bimodal distribution function. The data are divided
into two groups, each is described using a unimodal distribution
function, and a Beta pdf is utilized as \cite{mixed energy}:
 \vspace{-.5ex}
\begin{equation}
\label{eqn1} f_{b}(s)=\left\{
\begin{array}{l}
\frac{\Gamma(\alpha+\beta)}{\Gamma(\alpha)\Gamma(\beta)} \ast s^ {(\alpha+1)} \ast (1-s)^ {(\beta-1)}\\
0 \quad otherwise
\end{array}
\right.
\end{equation}
where $s$ is the solar irradiance ($kW/m^2$), $f_b(s)$ is the Beta
distribution function of $s$ and $\alpha$, $\beta$ are the parameters
of the Beta distribution function.

The output power of the PV module is dependent on the solar
irradiance and ambient temperature of the site as well as the
characteristics of the module itself. Therefore, once the Beta pdf per sitel
is generated for a specific time segment, the output power during
the different states is calculated for this segment using the
following:
 \vspace{-.5ex}
\begin{equation}
\label{eqn1}
\begin{array}{l}
T_{c_y}=T_A+s_{ay}\left( \frac{N_{OT}-20}{0.8}\right)\\
FF=\frac{V_{MPP}*I_{MPP}}{V_{oc}*I_{sc}}\\
I_y=s_{ay} \left[I_{sc} + K_i (T_{c_y} -25) \right]\\
V_y=V_{oc}-K_v*T_{c_y}\\
P_{s_y}(s_{ay})=N*FF*V_y*T_y
\end{array}
\end{equation}
where $T_{c_y}$ is the cell temperature $^\circ\mathrm{C}$ during
state $y$, $T_A$ is the ambient temperature $^\circ\mathrm{C}$,
$K_v$ is the voltage temperature coefficient V$/^\circ\mathrm{C}$;
$K_i$ is the current temperature coefficient A$/^\circ\mathrm{C}$;
$N_{OT}$ is the nominal operating temperature of cell in
$^\circ\mathrm{C}$; $FF$ is the fill factor; $I_{sc}$ is the short
circuit current in A; $V_{oc}$ is the open-circuit voltage in V;
$I_{MPP}$ is the current at maximum power point in A; $V_{MPP}$ is
the voltage at maximum power point in V; $P_{s_y}$ is the output
power of the PV module during state $y$; $s_{ay}$ is the average
solar irradiance of state $y$.

\subsubsection{Output Power of a Wind Turbine}

The wind speed profile is used to generate a corresponding output
power for use in the distributed power generation. Rayleigh pdf is
used to represent the wind speed profile for each time segment,
given as:
 \vspace{-.5ex}
\begin{equation}
\label{eqn2} f_{r}(v)=  \left(\frac{2v}{c^2}\right) \exp \left[
-\left(\frac{v}{c}\right)^ {2}\right]
\end{equation}
where $c=1.128v_m$.

The output power of a wind turbine is dependent on the wind speed at
the site as well as the parameters of its power performance curve.
Therefore, once the Rayleigh pdf is generated for a specific time
segment, the output power during the different states is calculated
for this segment as \cite{mixed energy}:
 \vspace{-.5ex}
\begin{equation}
\label{eqn1} P_{Vw}(v_{aw})=\left\{
\begin{split}
&0 & 0 \leq v_{aw} \leq v_{ci}\\
&P_{rated}*\frac{v_{aw}-v_{ci}}{v_{r}-v_{ci}} & v_{ci} \leq v_{aw} \leq v_{r}\\
&P_{rated} & v_{r} \leq v_{aw} \leq v_{co}\\
&0 & v_{co} \leq v_{aw}
\end{split}
\right.
\end{equation}
where $v_{ci}$, $v_r$ and $v_{co}$ are cut in speed, rated speed,
and cut-off speed of the wind turbine, respectively,
$P_{Vw}(v_{aw})$ is output power of the wind turbine during state
$w$, $v_{aw}$ is average wind speed of state.

\renewcommand{\thefootnote}{\fnsymbol{footnote}}
\footnotetext[1]{The value of all of the parameters related to PV
module and wind turbine is refer to \cite{mixed energy}.}

To incorporate the output power of the solar PV DG and wind turbine DG
units as multistate variables in the global formulation, the
continuous pdf of each has been divided into states (periods), in
each of which the solar irradiance and wind speed are within
specific limits. In other words, for each time segment, there are a
number of states for the solar irradiance and wind speed. The
probability of the solar irradiance and wind speed generation for
each state during any specific hour is calculated as follows:
 \vspace{-.5ex}
\begin{equation}
\label{eqn3} \Omega_{s}(G_y)= \int_{s_{y1}}^{s_{y2}} f_b(s)ds ,
\end{equation}
\begin{equation}
\label{eqn4} \Omega_{v}(G_w)= \int_{v_{w1}}^{v_{w2}} f_w(v)dv
\end{equation}
where $\Omega_{s}(G_y)$ is the probability of the solar irradiance
being in state $y$; $\Omega_{v}(G_w)$ probability of wind speed
being in state $w$; $s_{y1}$ and $s_{y2}$ are the solar irradiance
limits of state $y$; $v_{w1}$ and $v_{w2}$ are the wind speed limits
of state $w$.

Assuming that wind speed states and solar irradiance states are
independent, the probability of any combination of them
($\Omega\{C_g\}$) is obtained by convolving the two probabilities:
 \vspace{-.5ex}
\begin{equation}
\label{eqn5} \Omega\{C_g\}= \Omega_{s}(G_y)\ast \Omega_{v}(G_w).
\end{equation}

We define a supply function $S(t)$ to indicate the energy generated
by the renewable sources at time $t$, and obtain
 \vspace{-.5ex}
\begin{equation}
\label{eqn5} S(t)= \sum_{g=1}^G P_{g,t}\Omega\{C_g\}
\end{equation}
where $P_{g,t}$ is the mixed power of the $g$ state at time $t$. We
have a maximum limit on the electricity that can be generated by the
renewable plant i.e. $ 0 \leq S(t) \leq S_{max},\forall t.$

\subsection{Residential Smart Grid Model}

Let $\mathcal{N}$ denote the set of users in a residential network,
where the number of users, each representing one household
 is $N \triangleq \mid \mathcal{N} \mid$.
We consider the residential smart grid in which the electricity from
the power grid and renewable sources are shared by these $N$
users (homes), each one of which is equipped by a smart meter.
The smart meter connects with all the household appliances and not only
manages the electricity flow, but also schedules the power consumption of
each appliance in terms of the collected user's preference
information. Based on all the collected information, a centralized scheduler will
globally optimize the hourly power consumption and schedule all
appliances.

The appliances that define a household are divided into two types:
essential and shiftable appliances. Essential appliances are defined
as appliances that have a fixed energy requirement that would not be
subjected to scheduling decisions. These appliances include
television sets and refrigerators. The optimization process would
ensure a continuous supply of electricity for these. Shiftable
appliances allow the smart meter to schedule electricity loads to
minimize electricity cost. These appliances could include
dishwashers, washing machines and clothes dryers. We denote the
total essential and shiftable load for the $n$th user at time $t$ as
$a_{n,e}(t)$ and $a_{n,s}(t)$, respectively. Then, the total
essential and shiftable load across all users at time $t$ can be
respectively calculated as
 \vspace{-1ex}
\[A_e(t)=\sum_{n \in \mathcal{N}} a_{n,e}(t), \quad A_s(t)=\sum_{n \in \mathcal{N}}
a_{n,s}(t).\]

\subsection{Delay-Aware Power Transactions}
In our model, the energy produced by the distributed energy source at time $t$
can either be used to serve essential demand, shiftable demand or
can be sold to the grid. Let $S_e(t)$ and $S_s(t)$ denote the
electricity taken from the renewable plant to satisfy the essential
and shiftable demand at time $t$, respectively. Let $S_p(t)$ denote
the electricity that is sold by the customer to the grid at time
$t$. Then, we have
 \vspace{-1ex}
\[S_e(t)+S_s(t)+S_p(t)=S(t), \quad \forall t .\]
However, the renewable energy may not be enough to meet all of the
requests within a timely manner, and hence we also decide to
purchase an amount of electricity from the power grid. Let $G(t)$,
$G_e(t)$ and $G_s(t)$ denote respectively the total electricity
purchased from the grid, the portion of electricity taken from grid
to serve essential and shiftable loads at time $t$, we have,
 \vspace{-.5ex}
\[G(t)=G_e(t)+G_s(t), \quad \forall t .\]
Meanwhile, in order to formulate the problem mathematically, we put
a limit on max power that can be purchased from the grid, i.e. $ 0
\leq G(t) \leq G_{max},\forall t.$

Let $p(t)$ denote the cost of electricity taken from the grid at
time $t$. Let $\gamma(t)$ denote the price set by the grid to
purchase electricity at time $t$. The instantaneous cost of
electricity at time $t$, denoted by $C(t)$, is then given by,
 \vspace{-1ex}
\begin{equation}\label{inscost}
\begin{split}
C(t) &= p(t) G(t) - \gamma(t) S_p(t) \\
&= p(t)[G_e(t)+G_s(t)]  - \gamma(t) [S(t)-S_e(t)-S_s(t)].
\end{split}
\end{equation}

We assume that shiftable loads are flexible, and can tolerate their
demands being satisfied with some delay. We store all the shiftable
demands in a queue. The electricity demands in this queue are served
in a FIFO manner and the queue is then updated according to:
 \vspace{-1ex}
\[Q(t+1)=\max \bigg\{Q(t)-G_s(t)-S_s(t),0\bigg\}+A_s(t) .\]
Let $\bar{Q}$ denote the expected time average of the demand queue,
 \vspace{-1ex}
\[\bar{Q}=\text{lim}_{T\rightarrow \infty} \: \frac{1}{T} \sum_{t=0}^{T-1} E\bigg\{Q(t)\bigg\} < Q_{max}.\]The demands are served with finite delays if $Q_{max} <
\infty$.

\section{Problem Formulation}
Our objective now is to minimize the expected time average cost of
electricity subject to satisfying the essential demands and
guaranteeing the shiftable demands with a worst case delay bound
$D_{max}$. The optimization problem can now be formulated as
follows,
 \vspace{-1ex}
\begin{equation}\label{obj1}
P_1 = \min_{G_e(t),G_s(t),S_e(t),S_s(t)} \: \text{lim}_{T\rightarrow
\infty} \: \frac{1}{T}  \sum_{t=0}^{T-1}  E  \left\{ C(t) \right\}
\end{equation}
subject to,
\begin{equation}
G_e(t)+S_e(t)=A_e(t), \quad \forall t, \label{const1}
\end{equation}
\begin{equation}
0 \leq S(t) \leq S_{max}, \quad \forall t, \label{const3}
\end{equation}
\begin{equation}
0 \leq G(t) \leq G_{max}, \quad \forall t, \label{const4}
\end{equation}
\begin{equation}
Q(t+1)=\max \bigg\{Q(t)-G_s(t)-S_s(t),0\bigg\}+A_s(t),
\label{const5}
\end{equation}
\begin{equation}
\bar{Q} < Q_{max} \label{const6}
\end{equation}
where $C(t)$ is given by (\ref{inscost}). The backward
dynamic programming approach can be used to solve problem (\ref{obj1})
subject to the constraint set (\ref{const1})-(\ref{const6}). The state vector is given by
 $x_t=(S_t, A_{e,t}, A_{s,t}, p_t, \gamma_t)$,
where $S_t$ is the renewable supply, $A_{e,t}$ and $A_{s,t}$ are the
essential and shiftable demand respectively, and $p_t$, $\gamma_t$ are the market
price for buying and selling electricity respectively. Our control is assumed as
$u_t=(G_{e,t},G_{s,t},S_{e,t},S_{s,t})$, the backward dynamic
programming algorithm is given by the following equation
\cite{Neely1}:
\[J_t(x) = \min
_{u_t \in U_t(x)} \left\{C(x, u_t) + \sum_{ y \in
Y_{t+1}}J_{t+1}(y)\right\}\] where $J_t(x)$ is the value function of
period $t$, $C(x, u_t) = p_t(G_{e,t}+G_{s,t}) - \gamma_t
(S_t-S_{e,t}-S_{s,t})$ is the electricity cost incurred at each
period, $U_t(x)$ is the feasible region of actions for period $t$,
and $Y_t$ is the feasible region of the state vector at period $t$.

Since the dynamic programming need global information to solve the
problem which may have the disadvantage to capture the time-varying
characteristics of the demand and supply in renewable related
residential networks. It stands to reason that an algorithm that can
make optimal decision without the knowledge of future information
should be designed. Inspired by \cite{Neely2}, we design a
distributed online algorithm which will utilize a Lyapunov
optimization technique to solve the optimization problem in a
distributed way without a-priori statistical knowledge.

\section{Distributed Online Algorithm}

\subsection{Lyapunov Optimization}

Let $P_1^*$ denote the infimum time average cost of the above
problem. The requirement of $\bar{Q} < \infty$ in the above problem
does not guarantee a delay bound $D_{max}$. To ensure this bound a
delay aware virtual queue is developed. This delay aware virtual
queue has the same service rate as that of the actual queue but it
grows if the demand is not serviced at time $t$. The virtual queue
update equation is,
 \vspace{-.5ex}
\begin{equation}
Z(t+1)=\max \bigg\{Z(t)-G_s(t)-S_s(t)+\epsilon 1_{\big\{Q(t)
>0\big\}} ,0 \bigg\} \label{const7}
\end{equation}
where, $1_{\big\{Q(t) >0\big\}}$ is an indicator function. The value
of this function is 1 if $Q(t) >0$ and 0 otherwise. Thus the delay
aware queue grows by an amount $\epsilon$ each time it is not
serviced. With the introduction of delay aware virtual queue we have
the following Lemma \cite{Neely1}. \\
 \vspace{-1ex}
\\Lemma 1: If $Z(t) \leq Z_{max}$ and $Q(t)\leq Q_{max}$, then all shiftable devices will get scheduled
with a worst case delay of $D_{max}$ where,
 \vspace{-.5ex}
\[D_{max}=\left\lceil \frac{Z_{max}+Q_{max}}{\epsilon} \right\rceil .\]
Let $\Theta(t)=(Z(t),Q(t))$ denote the concatenated vector of
virtual and real queues. We then define the following Lyapunov
function to measure the congestion in the queues,
 \vspace{-.5ex}
\[L(\Theta(t))=\frac{1}{2} [Z^2(t)+Q^2(t)] .\]
The conditional 1-slot Lyapunov drift is defined as,
 \vspace{-.5ex}
\[\Delta(\Theta(t))=E\bigg\{L(\Theta(t+1))-L(\Theta(t)) | \Theta(t) \bigg\} .\]
Define $V>0$ as a parameter to effect the performance delay
tradeoff, we have the following drift plus penalty function, by
using the drift plus penalty framework developed in \cite{Neely1},
 \vspace{-.5ex}
\[\Delta(\Theta(t)) + V E\bigg\{ C(t) | \Theta(t)  \bigg\} \]
where $C(t)$ is given by (\ref{inscost}). Again using the same
techniques in \cite{Neely1} we can show that this function is
bounded as,
 \vspace{-.5ex}
\[\Delta(\Theta(t)) + V E\bigg\{ C(t) | \Theta(t)  \bigg\} \leq B+ V E\bigg\{ C(t) | \Theta(t) \bigg\} \]
\[+Q(t) E\bigg\{ \big(A_s(t)-G_s(t)-S_s(t) \big) | \Theta(t) \bigg\} \]
\begin{equation}\label{bound}
+Z(t) E\bigg\{ \big(\epsilon-G_s(t)-S_s(t) \big) | \Theta(t)
\bigg\}.
\end{equation}
Now instead of solving the original problem we minimize this bound
on the drift plus penalty function.

\subsection{Two-Stage Online Algorithm}

In this section, we propose a two-stage online algorithm to minimize
the right-hand-side of the drift-plus-penalty bound (\ref{bound}) at
every slot $t$ as follows,\\
\\1. At each time $t$ observe: $Z(t)$, $Q(t)$, $p(t)$, $\gamma(t)$, $S(t)$, $A_s(t)$ and $A_e(t)$.
\\2. Then solve the following optimization problem to obtain the decision variables $\{G_e^*(t),G_s^*(t),S_e^*(t),S_s^*(t) \}$:
\begin{eqnarray}
\min_{G_e(t),G_s(t),S_e(t),S_s(t)} \: V C(t)  \nonumber
-\bigg\{Q(t)+Z(t)\bigg\} \bigg\{G_s(t)+S_s(t) \bigg\} \label{obj2}
\end{eqnarray}
subject to constraints (\ref{const1}), (\ref{const3}),
(\ref{const4}).
\\3. Then update the queue equations for next time slot,
 \vspace{-.5ex}
\[Q(t+1)=\max \bigg\{Q(t)-G_s^*(t)-S_s^*(t),0\bigg\}+A_s(t), \]
 \vspace{-1ex}
\[Z(t+1)=\max \bigg\{Z(t)-G_s^*(t)-S_s^*(t)+\epsilon 1_{\big\{Q(t) >0\big\}} ,0 \bigg\} .\]

Using this algorithm, the problem is decoupled in time. Also we do
not require any information about the future values or the
expectations. In addition, we observe that the electricity served
for essential loads $G_s(t)$ and $S_e(t)$ are independent with that
for shiftable load in the above formulation. Hence, our problem can
be separately solved in two stages. In stage I, we will minimize the
cost and only introducing essential load as the input; In stage II,
we solve the optimization problem for shiftable load.

\subsubsection{Stage I, Optimization for Essential Load}

At each time slot the instantaneous optimization problem is,
 \vspace{-1ex}
 \begin{eqnarray}
\min_{G_e(t),\!G_s(t),\!S_e(t),\!S_s(t)} \: \! \! V C(t) \nonumber
 \!-\!\bigg\{Q(t)\!+\!Z(t)\bigg\} \!\bigg\{G_s(t)\!+\!S_s(t) \!\bigg\} \label{obj3}
\end{eqnarray}
subject to constraints (\ref{const1}), (\ref{const3}),
(\ref{const4}). Rearranging the terms in the above objective
function we get the following optimization problem for the essential
load,
 \vspace{-.5ex}
\begin{equation}
\min_{G_e(t),S_e(t)} \: V\bigg(p(t) G_e(t) + \gamma(t) S_e(t) \bigg)
\label{ess_obj}
\end{equation}
subject to,
 \vspace{-.3ex}
\begin{equation}
G_e(t)+S_e(t)=A_e(t), \quad \forall t. \label{essconst1}
\end{equation}
Replace $G_e(t)=A_e(t)-S_e(t)$ in the above objective function, we
get the following unconstrained optimization problem in only one
variable,
 \vspace{-.5ex}
\begin{equation}
\min_{S_e(t)} \: V(\gamma(t)-p(t)) S_e(t). \label{essobj2}
\end{equation}
From this equation it is obvious that for cost minimization we
should adopt the following strategy,
 \vspace{-1ex}
\begin{eqnarray}
G_e^*(t)=\Bigg\{\begin{array}{c}
                \tilde{G}=\min \big(G_{max}, A_e(t)\big),\quad \text{if} \: p(t) < \gamma(t) \\
                \max \big(A_e(t)-\tilde{S},0\big),\quad \text{if} \: p(t) \geq
                \gamma(t).
                \end{array}
\label{ess_grid}
\end{eqnarray}
 \vspace{-1ex}
\begin{eqnarray}
S_e^*(t)=\Bigg\{\begin{array}{c}
                \tilde{S}=\min \big(S(t), A_e(t)\big),\quad \text{if} \: \gamma(t) \leq p(t) \\
                \max \big(A_e(t)-\tilde{G},0\big),\quad \text{if} \: \gamma(t) >
                p(t).
                \end{array}
\label{ess_ren}
\end{eqnarray}
This result is very simple and it says that we should use maximum
power from the grid if the price of purchasing electricity is less
than price of selling renewable electricity back to the grid.
Otherwise utilize maximum power from the renewable source to satisfy
the essential load demands. We assume that there is always enough
power available to serve the essential demands in any time slot.

\subsubsection{Stage II, Optimization for Shiftable Load} Once the
decision on $G_e^*(t)$ and $S_e^*(t)$ for essential load are
obtained we have to decide about
the remaining renewable and grid power. The power that is available for shiftable demands is, 
 \vspace{-.5ex}
 \begin{equation}
\tilde{S}_{max}=\max(S(t)-S_e^*(t),0), \label{const13}
\end{equation}
 \vspace{-4ex}
\begin{equation}
\tilde{G}_{max}= \max(G_{max}-G_e^*(t),0). \label{const14}
\end{equation}
The objective function for the shiftable problem now contains the
remaining terms and the remaining constraints. After some
re-arrangement of terms in the objective function the instantaneous
optimization problem to obtain variables $G_s(t)$ and $S_s(t)$ is
now as follows,
 \vspace{-1ex}
\begin{eqnarray}
\min_{G_s(t),S_s(t)} \: \bigg\{ G_s(t) \bigg(V
p(t)-(Q(t)+Z(t))\bigg) \bigg\} + \nonumber \\ \bigg\{ S_s(t) \bigg(V
\gamma(t)- (Q(t)+Z(t)) \bigg) \bigg\} \label{inst_obj}
\end{eqnarray}
subject to constraints  (\ref{const13}), (\ref{const14}). We can
then decide the optimum variables based on appropriate threshold
type policy.

\subsection{Real-Time Distributed Decision}
There are four possible cases depending on the current prices and
queue states.
\subsubsection{Case I: $V p(t) > Q(t)+Z(t)$ and $V \gamma(t) > Q(t)+Z(t)$}
This situation indicates any or all of the following:
\begin{itemize}
    \item the queue lengths are small
    \item for given $V$, $p(t)$ is large i.e. price of buying electricity from grid is high
    \item for given $V$, $\gamma(t)$ is large i.e. price of selling electricity to grid is high
\end{itemize}
In this situation (also evident from the objective function),
shiftable load should not be scheduled and instead we should sell
all the electricity to the grid i.e.
\[G_s^*(t)=0, \quad S_s^*(t)=0,\quad S_p^*(t)=\tilde{S}_{max}.\]

\subsubsection{Case II: $V p(t) > Q(t)+Z(t)$ and $V \gamma(t) < Q(t)+Z(t)$}
This situation is possible when $p(t) \geq \gamma(t)$. In this case,
we should serve the shiftable queue using renewable energy only. The
shiftable queue should not be served from the grid and any remaining
renewable energy left after servicing the shiftable queue should be
sold to the grid.
 \vspace{-1ex}
\[G_s^*(t)=0,  \]  \vspace{-4ex}
\[S_s^*(t)=\min (\tilde{S}_{max},Q(t)), \] \vspace{-4ex}
\[S_p^*(t)=\max (\tilde{S}_{max}-S_s^*(t),0). \]

\subsubsection{Case III: $V p(t) < Q(t)+Z(t)$ and $V
\gamma(t) > Q(t)+Z(t)$} This situation is possible when $p(t) <
\gamma(t)$. This case results in a situation where the customer
would sell all the renewable electricity to the grid and purchase
electricity from the grid to serve the shiftable demands in the
queue. However, to prevent users to buy energy from the grid and
sell it back, a rational pricing model will make sure that $p(t)$
will never less than $\gamma(t)$. Hence, this case will not happen.


\subsubsection{Case IV: $V p(t) < Q(t)+Z(t)$ and $V \gamma(t) < Q(t)+Z(t)$}
This situation means any or all of the following:
\begin{itemize}
    \item queue thresholds are high
    \item price of electricity from grid is low
    \item price offered by grid for renewable energy is also low
\end{itemize}
In this case we should serve the queue using the energy from the
grid as well as the renewable source. For this situation we have two
further sub-cases: \\
\emph{Sub-case I} (if $p(t) \geq \gamma(t)$): in this case, the
consumers use maximum energy from the renewable source to serve
shiftable loads. The remaining demands in the queue should then be
served using grid energy. If some renewable energy remains after
servicing the queue then it should be sold to the grid.
 \vspace{-1ex}
\[S_s^*(t)=\min (\tilde{S}_{max},Q(t)), \] \vspace{-3ex}
\[G_s^*(t)=\min (Q(t)-S_s^*(t),\tilde{G}_{max}),  \]  \vspace{-3ex}
\[S_p^*(t)=\max (\tilde{S}_{max}-S_s^*(t),0). \]
\emph{Sub-case II} (if $p(t) < \gamma(t)$): Based on the same reason
expressed in Case III, this case also will not happen.

\section{The Simulation Analysis}

We consider 10 customers (homes) and each customer is selected to
have 20 essential appliances and 20 shiftable appliances. We have a
set up similar to \cite{Demand side}, where the essential appliances
may include electric stoves (daily usage $E$: 1.89 kWh for
self-cleaning and 2.01 kWh for regular), lighting (daily usage for
10 standard bulbs: 1.00 kWh). Based on this set up, the essential
load $A_e(t)$ for each hour can be obtained. The shiftable
appliances may include dishwashers (daily usage $E$: 1.44 kWh),
washing machines (daily usage $E$: 1.49 kWh for energy-start, 1.94
kWh for regular) and clothes dryers (daily usage $E$: 2.50 kWh). To
describe the time-varying characteristic of the shiftable load, we
model the shiftable demand $A_s(t)$ as i.i.d. over slots and
uniformly distributed over the integers $\{0, 1, \cdots,
A_{s,max}\}$. We assume that market prices for purchasing
electricity from the grid is $p(t)=0.3$ cents at daytime hours,
i.e., from 8:00 in the morning to 12:00 at night and $p(t)=0.1$
cents during the night, i.e., from 12:00 at night to 8:00 AM the day
after. Meanwhile, the selling price is set at $\gamma(t)=0.2$ cents
during daytime hours and $\gamma(t)=0.1$ cents during the night.

\subsection{Electricity Cost}
We compare the proposed bidirectional transaction scheme with
dynamic programming (BTS-DP) and Lyapunov optimization (BTS-LO)
against a simple ``Purchasing Only Strategy (POS)'' which tries to
use all the available renewable resource $S(t)$ and only buy from
the market as a last resort if the produced renewable energy cannot
satisfy the current essential and shiftable load. Fig. \ref{fig:
ec1} shows the cumulative electricity costs in different power
transaction schemes. It is shown that the energy costs achieved by
the proposed bidirectional transaction schemes are less than that in
the POS. This is because the proposed schemes are not only able to
sell the produced renewable energy to the grid to obtain financial
credits, but also schedule the shiftable loads to the lower
electricity price periods to reduce electricity consumption costs.
We note that the total electricity cost in the BTS-LO case is higher
than that in the BTS-DP. That is, the BTS-LO is a distributed online
optimization and requires less information than BTS-DP. However, we
can obtain the different costs in BTS-LO by adjusting the parameters
$V$ and $\epsilon$, which is related to the delay suffered by
shiftable load.

 \vspace{-3ex}
\begin{figure}[htb]
\centerline{\includegraphics[width=7.5cm]{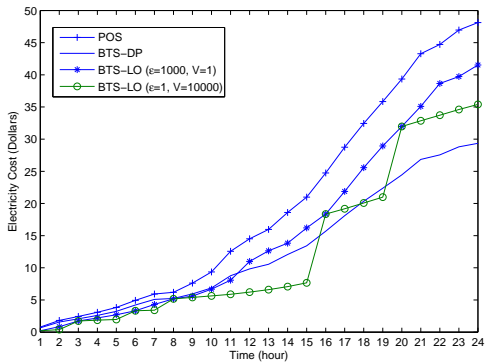}}
\vspace{-3ex} \caption{Electricity cost in terms of different power
transaction schemes.} \label{fig: ec1}
\end{figure}

 \vspace{-1ex}
\subsection{Delay}
\begin {figure}[htb]
\centerline{\includegraphics[width=7.5cm]{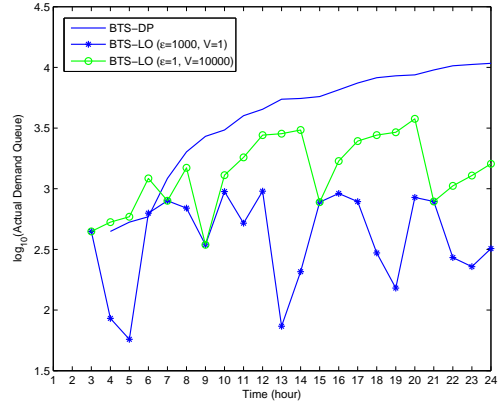}}
\vspace{-3ex}\caption{Delay comparison between the bidirectional
transaction scheme with dynamic programming and Lyapunov
optimization.} \label{fig: delay1}
\end{figure}

To illustrate the delay in the proposed scheme, we use $log_{10}$
(Actual Demand Queue ($Q$)) to compare the delay between BTS-DP with
infinite delay constraint and BTS-LO with delay constraint
($D_{max}$). The delay in terms of different scenarios have shown in
Fig. \ref{fig: delay1}. In Fig. \ref{fig: delay1}, it is note that
the $log_{10}$ (Actual Demand Queue ($Q$)) has the highest value in
the dynamic programming solution, as it achieves the lowest energy
cost. In addition, selecting different parameters $V$ and $\epsilon$
in real-time algoirthm BTS-LO, we can adjust the delay and achieve
corresponding electricity costs (as shown in Fig. \ref{fig: ec1}),
and thus, achieve various tradeoff of the consumers' electricity
usage preference: power consumption vs delay tolerance.

 \vspace{-1ex}

\section{Conclusion}
In this paper, we presented a bidirectional power transaction
mechanism to minimize the electricity cost in the residential smart
grid. According to the endurance of the serving delay, the
residential load is classified into two types: essential load and
shiftable load. In order to satisfy both of the loads, consumers can
purchase the power from grid or use the produced renewable power.
Since the consumers also can sell the produced renewable power to
grid for benefit, the shiftable load can be stacked in a queue,
given that it will be served within a pre-agreed deadline. To decide
the amount of electricity consumed, sold and purchased, we then
formulate a global centralized optimization problem and provide a
distributed online algorithm by using Lyapunov optimization
solution. Simulations confirms that our proposed mechanism is able
to balance the tradeoff between the power consumption and the
allowable delay in the residential networks. Meanwhile, the
effectiveness of distributed online algorithms with limited
information is compared with the bound provided by the centralized
off-line algorithm with global information.

 \vspace{-.5ex}
\section{Acknowledgment}
This research is partly supported by SUTD-ZJU/RES/02/2011,
International Design Center, EIRP, and SSE-LUMS via faculty research
startup grant.

 \vspace{-.5ex}
\appendices

\newpage

\end{document}